\documentclass[12pt]{article} 

\usepackage{amsmath}
\usepackage{amsthm}
\usepackage{amsfonts}
\usepackage{mathrsfs}
\usepackage{stmaryrd}
\usepackage{setspace}
\usepackage{fullpage}
\usepackage{amssymb}
\usepackage{breqn}
\usepackage{enumitem}
\usepackage{bbold} 
\usepackage{authblk}
\usepackage{comment}
\usepackage{hyperref}
\usepackage{pgf,tikz}
\usepackage{graphicx}
\usepackage{subcaption}

\newtheorem*{thm*}{Theorem}
\newtheorem{thm}{Theorem}[section]

\newtheorem{proposition}[thm]{Proposition}
\newtheorem*{prop*}{Proposition}

\newcommand\cF{{\mathcal F}}

\newcommand\cH{{\mathcal H}}

\newcommand{\ignore}[1]{}

\title{Cooperation in Combinatorial Search
}
\author[1]{D\'aniel Gerbner} \author[1,2]{Bal\'azs Keszegh} 
\author[2]{Kartal Nagy} 
\author[1]{Bal\'azs Patk\'os} \author[3]{G\'abor Wiener}
\affil[1]{Alfr\'ed R\'enyi Institute of Mathematics}
\affil[2]{ELTE Eötvös Lor\'and University}
\affil[3]{Budapest University of Technology and Economics,  Department of Computer Science and Information Theory}

\date{}

\begin{document}

\maketitle

\begin{abstract}
    In the game theoretical approach of the basic problem in Combinatorial Search an adversary thinks of a defective element $d$ of an $n$-element pool $X$, and the questioner needs to find $x$ by asking questions of type is $d\in Q?$ for certain subsets $Q$ of $X$. We study cooperative versions of this problem, where there are multiple questioners, but not all of them learn the answer to the queries. We consider various models that differ in  how it is decided who gets to ask the next query, who obtains the answer to the query, and who needs to know the defective element by the end of the process. 
\end{abstract}

\section{Introduction}

In Combinatorial Search Theory, the basic problem is the following. We are given an $n$-element underlying set $X$ that contains an unknown defective element $d$ (we assume $n\ge 2$, otherwise all search problems are straightforward). We have to identify the defective element by asking queries of type is $d\in Q?$ for certain subsets $Q \subseteq X$ (we refer to such a query simply as asking $Q$). Our aim is of course to use as few queries as possible in order to identify $d$. What is the minimum number of questions needed and how can we find the queries corresponding to this minimum?

This basic question is easy to answer if we can ask any subset of $X$: we need  $\lceil \log n\rceil$ queries (here and throughout the paper $\log n$ denotes $\log_2 n$). There are several variants that lead to interesting questions: there might be a restriction on the query sets, some answers can be erroneous, the queries have to be asked at the same time or in a certain number of batches (non-adaptive and $k$-round algorithms), there might be more defectives, and so on; see \cite{DH1999} for more details.

Here we introduce a new family of variants, which involve cooperation. Rather than one participant asking the queries and trying to identify the defective element, we are dealing with more participants that have to work together. We consider several models differing in mainly the information the participants obtain after a query. We use the game theoretic approach of Search Theory (see \cite{A}), where the search process is considered to be a game between Player $A$ (the adversary) and Player $B$ (the questioner) and extend it by adding more players.

Our starting point is a nice puzzle by S\'andor R\'oka \cite{bhr} (motivated by his research in \cite{RR,RR2}). Player $A$ picks the defective element and answers the queries, Player $B$ asks the queries (all at the same time, non-adaptively) and Player $C$ is given the YES answers (along with the sets to which these answers correspond, of course). In other words, Player $B$ chooses a family $\cF$ of subsets of $X$, and Player $A$ gives the subfamily $\cF'$ to Player $C$, where $\cF'=\{F\in \cF: \, d\in F\}$. Upon receiving the answers, Player $C$ has to identify $d$, which can be done if and only if there exists a unique element of $X$ contained in each member of $\cF'$.

It is not hard to see that Player $C$ can identify $d$ if and only if $\cF$ is \textit{completely separating}, i.e., $\cF$ satisfies the following property: for any $x,y\in X$, there exists a set $F\in \cF$, such that $x\in F$, $y\not\in F$. The smallest possible size of a completely separating family was determined by Spencer \cite{spen}. As we will use its ideas later, we briefly describe the simple proof.

The \textit{dual} of a family $\cF$ of subsets of $X$ is defined as a family $\cH$ of cardinality $|X|$ on the underlying set $\cF$ as follows. For each $x\in X$, we define a set $H_x\in \cH$ such that for each $F\in \cF$ we have $F\in H_x$ if and only if $x\in F$. It is easy to see that $\cF$ is completely separating if and only if its dual is a \textit{Sperner system}, i.e., there are no $H_x$ and $H_y$ in $\cH$, such that $H_x\subset H_y$. The classical theorem of Sperner \cite{S} states that the largest possible cardinality of a Sperner system on an $m$-element underlying set is $\binom{m}{\lfloor m/2\rfloor}$. This implies that the smallest possible size of a completely separating family on an $n$-element underlying set is the smallest $m$.such that $\binom{m}{\lfloor m/2\rfloor}\ge n$.

Let us now turn to new models. In each of the models, Player $A$ picks the defective element and answers the queries. Some models will require adaptive, others require non-adaptive algorithms. Models also differ by numerous other properties, like 
\begin{itemize}
\item whether the players can agree in a strategy beforehand, 
\item whether just one player can ask queries or more, and in the latter case what is the order in which the players are allowed to ask their queries, 
\item what information the players obtain after a query,
\item which players should be able to identify the defective element.
\end{itemize}

We assign numbers to the models and say that an algorithm solves Model $i$ if the required players can always identify the defective element. If such an algorithm  exists, the minimum number of queries needed for the solution (in the worst case, of course) is denoted by $f_i(n)$. Since all models we consider are harder to solve than the basic model of search, the trivial lower bound $f_i(n)\ge \lceil \log n\rceil$ holds for all $i$.

We start with some models where Player $B$ asks the queries and multiple players obtain different pieces of information from the answers. These models are non-adaptive, and we assume that the players other than Player $B$ do not have any information beforehand about the algorithm: they are only given a family of sets with the information whether they contain the defective element.

For each model we state our results and prove these later in the next section.

\bigskip

\noindent\textbf{Model 1.}

We have four players $A,B,C,D$. Player $B$ asks the queries non-adaptively, Player $C$ obtains the YES answers and Player $D$ obtains the NO answers (together with the corresponding queries). At least one of $C$ and $D$ have to be able to identify the defective element. 

In order to find the optimum strategy for Model 1, we use the notion of $k$-Sperner systems. A family is said to be $k$-Sperner if it does not contain $k+1$ distinct sets $S_1,\dots,S_{k+1}$ such that $S_1\subsetneq S_2\subsetneq\dots\subsetneq S_{k+1}$. A theorem of Erd\H os \cite{Er} states that the largest possible cardinality of a $k$-Sperner family on an $m$-element underlying set is $\sum_{i=1}^k\binom{m}{\lfloor \frac{m-k}{2}\rfloor+i}$, from which the second statement of the following Proposition follows.

\begin{proposition}\label{prop1} A family $\cF$ solves Model 1 if and only if the dual $\cH$ of $\cF$ is a 2-Sperner family. Therefore, $f_1(n)=\min \{m: \binom{m}{\lfloor m/2\rfloor}+\binom{m}{\lfloor m/2+1\rfloor}\ge n\}$.
\end{proposition}

\bigskip

Our next model is a generalization of Model 1.

\smallskip

\noindent\textbf{Model 2.}

We have $r+2$ players, $A,B,C_1,\dots,C_r$ with $r>1$. Player $B$ asks a partition of $X$ to $r$ sets $X_1,\dots,X_r$ as a query. The defective element is contained in $X_i$ for some $i$; then Player $C_i$ is given a YES answer (together with what the query was). We also say that $X_i$ belongs to Player $C_i$. The goal is that  at least one of Players $C_j$ should be able to identify the defective element. Note that Model 2 is a generalization of Model 1, indeed, since a query $F$ of $B$ in Model 1 corresponds to the query $(F,X\setminus F)$ of $B$ in Model 2.

It is more convenient to state our result on $f_2(n)$ (or,  more precisely, on $f_2(n,r)$) in the following form. Let $g_2(k,r)$ denote the largest number of elements such that $k$ queries are enough to solve Model 2 on $g_2(k,r)$ elements. Then clearly (as in the previous models) $f_2(n,r)$ is the smallest $k$ such that $n\le g_2(k,r)$.

\begin{thm}\label{thm2}
    \begin{displaymath}
g_2(k,r)=
\left\{ \begin{array}{l l}
r\binom{k}{\lfloor k/2\rfloor} & \textrm{if\/ $k\ge 3$ and $r\neq 2$},\\
\binom{k}{\lfloor k/2\rfloor}+\binom{k}{\lfloor k/2+1\rfloor} & \textrm{if\/ $r=2$}.\\
2r-1 & \textrm{if\/ $k=2$}.\\
\end{array}
\right.
\end{displaymath}
\end{thm}


\bigskip

The next model differs from Model 2 only in that Player $B$ can ask $r$ disjoint sets that do not need to form a partition. This also means that $r=1$ makes sense now and in this case we obtain the original puzzle we have dealt with earlier.

\medskip

\noindent\textbf{Model 3.}

We are given $r+2$ players, $A,B,C_1,\dots,C_r$. Player $B$ asks $r$ disjoint subsets $X_1,\dots,X_r$ of $X$. If the defective element is contained in $X_i$ for some $i$, then Player $C_i$ obtains a YES answer (together with what the query was). The goal is that  at least one of Players $C_i$ should be able to identify the defective element.

\begin{thm}\label{thm3}
    $g_3(k,r)=r\binom{k}{\lfloor k/2\rfloor}$.
\end{thm}

Let us turn our attention to some completely different models. From now on we will have multiple questioners and the goal is to minimize the total number of queries they ask. This makes sense only if they can affect each other, thus here we deal with adaptive search. We assume in the sequel that Players $B$ and $C$ can discuss their (questioning) strategies before asking any queries (but any communication between them is forbidden during the algorithm, of course). We consider only YES-NO queries, where only the YES answer is shared with some player. We also assume that this player does not know whether a query with a NO answer was asked (otherwise they would know that the answer is NO just from hearing no answer). Since such a player may need to ask the next query, the order in which the players ask the queries should be hidden from some players. There is an easy way to achieve this: Player $A$ decides which of the other players will ask the next query.

\bigskip
\noindent\textbf{Model 4.}

We have three players $A,B,C$. Player $A$ decides who asks next. Whenever Player $B$ asks a query, if the answer is YES, Player $C$ obtains what the query was (and that the answer was YES). Similarly, whenever $C$ asks a query, if the answer is YES, $B$ obtains what the query was  (and that the answer was YES). Note that the player asking the query does not get any information concerning this particular query. The goal is that at least one of $B$ and $C$ should be able to identify the defective element.

\begin{thm}\label{thm4}
    $f_4(n)\le \lceil\log n\rceil+\lceil\log\log n\rceil+2$.
\end{thm}

\bigskip
\noindent\textbf{Model 5.}

We have three players $A,B,C$. Player $A$ decides who asks next. Whenever Player $B$ asks a query, $B$ obtains the answer and if the answer is YES, $C$ also obtains what the query was (and that the answer was YES). Similarly, whenever $C$ asks a query, $C$ obtains the answer and if the answer is YES, Player $B$ also obtains what the query was (and that the answer was YES). The goal is that at least one of $B$ and $C$ should be able to identify the defective element.
Note that Model 5 differs from Model 4 in that the player asking the query always obtains the answer.

\begin{proposition}\label{prop5}
$f_5(n)\le \lceil\log n\rceil+1$.
\end{proposition}

This upper bound is only one larger than the trivial lower bound. For small values of $n$, we determined the exact value of $f_5(n)$ and it turns out that both $\lceil\log n\rceil$ and $\lceil\log n\rceil+1$ might occur as such, but we were not able to recognize a pattern. It seems to be an interesting question to determine the set of those $n$'s for which $f_5(n)=\lceil\log n\rceil$.

\bigskip
\noindent\textbf{Model 6.}

This model is the same as Model 5 with only one difference: the goal is that both $B$ and $C$ should be able to identify the defective element.

\begin{thm}\label{thm6}
    $f_6(n)\le \lceil\log n\rceil+2\lceil\sqrt{\log n}\rceil+2$.
\end{thm}

Finally, we introduce a third type of models. In these we handle a problem we have mentioned earlier in a different way: if a player obtains information about a query only if the answer is YES, then in the case of a NO answer they should not know what the query was. 
Now we handle this problem by assuming that the elements are indistinguishable before the start of the algorithm. That is, the players can decide that the first query set (say) is of order $n/2$, but cannot decide which $n/2$ elements should be included in the query set. Then they can decide that the next query contains (say) $n/4$ elements from the first query set and $n/3$ from the complement of the first query set, but they have no way to describe which elements are in the queries. This way, even if Player $B$ knows that the answer was NO to a query asked by Player $C$, this information is much less useful as before.

\bigskip
\noindent\textbf{Model 7.}

The elements are indistinguishable. Player $B$ and Player $C$ alternate asking queries. From here we follow Model 4: Whenever Player $B$ asks a query, if the answer is YES, Player $C$ obtains what the query was. Similarly, whenever Player $C$ asks a query, if the answer is YES, Player $B$ obtains what the query was. Note that the player asking the query does not obtain any information. The goal is that at least one of Player $B$ and Player $C$ should be able to identify the defective element. We remark that if the answer is NO to a query asked by (say)  $B$, then $C$ finds out that the answer was NO, since he has to ask the next query without obtaining a YES answer. However, Player $C$ does not find out what the query of Player $B$ was (this is why in this model we need indistinguishable elements, otherwise the two players could emulate the classic non-adaptive  search strategy using at most $\lceil\log n\rceil+1$ steps).


\bigskip
\noindent\textbf{Model 8.}

Model 8 is the same as Model 7, except that the player asking a query also obtains the answer to the query. Clearly $f_8(n)\le f_7(n)$.

Recall that the golden ratio is $\phi=\frac{1+\sqrt{5}}{2}$. 

\begin{thm}\label{thm7}
    $\log_\phi n -O(1)\le f_8(n)\le f_7(n)\le\lceil\log_\phi(n)\rceil$.
\end{thm}

\section{Proofs}

Let us start with the proof of Proposition \ref{prop1}. Recall that it states that a family $\cF$ solves Model 1 if and only if its dual $\cH$ is 2-Sperner.

\begin{proof}[Proof of Proposition \ref{prop1}]
Let us assume that $\cF$ solves Model 1. If Player $C$ can find the defective element $d$, then for every element $x\neq d$, there is a set $F\in \cF$ with $d\in F$, $x\not\in F$. If Player $D$ can find the defective element, then for every element $x\neq d$ there is a set $F'$ with $x\in F'$, $d\not\in F'$. As any element $y$ can be defective, we have for every $y$ that either for every $x\neq y$ there is a set $F$ with $y\in F$, $x\not\in F$, or for every $x\neq y$ there is a set $F'$ with $y\not\in F'$, $x\in F'$. 

It is equivalent that in the dual family $\cH$ for every set $H_y$, we have that either every other set $H_x$ has the property that there is an $F$ with $F\in H_y$, $F\not\in H_x$, or every set $H_x$ has the property that there is an $F'$ with $F'\not\in H_y$, $F'\in H_x$. The first property means that $H_y$ is not contained in any $H_x$, and the second property means that $H_y$ does not contain any $H_x$. In other words, $H_y$ is either maximal or minimal with respect to containment. Three distinct sets with $H\subset H'\subset H''$ in $\cH$ would mean that $H'$ violates this property, thus $\cH$ is 2-Sperner.

Observe that in each step of the above proof the implications go both ways, giving us the other direction. Finally, the result on $f_1(n)$ follows from the theorem of Erd\H os on $k$-Sperner families, mentioned in the introduction.
\end{proof}

We prove Theorems \ref{thm2} and \ref{thm3} together, as parts of their proofs are common. Observe that by definition we have $g_3(k)\ge g_2(k)$.

\begin{proof}[Proof of Theorems \ref{thm2} and \ref{thm3}]
    Let us start with proving the general upper bound $r\binom{k}{\lfloor k/2\rfloor}$. Observe that if $\cF$ is a solution to Model 2 or 3, then for every $x\in X$ we have at least one player $C_i$ such that if $x$ is the defective element, then $C_i$ can identify it. Let $Y_i$ be the set of elements that can be identified by $C_i$. 
    
    We claim that $|Y_i|\le \binom{k}{\lfloor k/2\rfloor}$. Let $x,y\in Y_i$, then there is a query where $x\in X_i,y\not\in X_i$, otherwise if $x$ is the defective, then both $x$ and $y$ appear in all the sets $X_i$ where Player $C_i$ obtains any information, thus Player $C_i$ does not know whether $x$ or $y$ is the defective element. This means that the sets $X_i$ restricted to $Y_i$ form a completely separating family $\cF_i$. Therefore, the dual family of $\cF_i$ is a Sperner family of cardinality $|Y_i|$ on the underlying set $\cF_i$ with $|\cF_i|=k$, completing the proof of the claim.

    Summing this we get that the number of elements is indeed at most $\sum |Y_i|\le r\binom{k}{\lfloor k/2\rfloor}$.

    Let us continue with the proof of the upper bound in Theorem \ref{thm2} in the case $k=2$. Assume on the contrary that $g_2(2,r)=2r$, in this case by the above claim we must have $|Y_i|=2$ for each $i$ and these sets are pairwise disjoint. Let $X_1,\dots,X_r$ be the first query and $X'_1,\dots,X'_r$ be the second query. Let $Y_1=\{x,y\}$, then clearly $x\in X_1,y\not\in X_1$ and $y\in X'_1,x\not\in X'_1$ (or the other way around) because of the completely separating property. As we deal with Model 2, another set in the first query, say $X_2$, contains $y$. Let $Y_2=\{u,v\}$ with $u\in X_2$. Then $u\ne y$ as $Y_1\cap Y_2=\emptyset$. If $u$ is the defective, then Player $C_2$ cannot rule out the possibility that $y$ is the defective, since Player $C_2$ only obtains the information that $X_2$ contains the defective (recall that in Model 2 and 3 the Players $C_i$ are not aware of the queries of Player $B$, they only get some sets that contain the defective element), a contradiction. This completes the proof of the upper bound $g_2(2,r)\le 2r-1$.

    Let us continue with the lower bound $r\binom{k}{\lfloor k/2\rfloor}$ in the case of Theorem \ref{thm3}. Consider a set $X$ of size $r\binom{k}{\lfloor k/2\rfloor}$. Player $B$ first partitions $X$ to $r$ sets $Y_1,\dots,Y_r$, each of size $\binom{k}{\lfloor k/2\rfloor}$. On each $Y_i$, Player $B$ picks a completely separating family of order $k$ and let $X_1^{i},\dots,X_k^{i}$ be its members. Then as Query $Q_j$, Player $B$ asks the sets $X_j^{i}$. If the defective element $d$ belongs to $Y_i$, then $C_i$ can identify $d$, as any element not in $Y_i$ is distinguished by any $X^i_j$ that contains $d$ (and there is at least one such set) and any $y\neq d$ in $Y_i$ is distinguished from $d$ as $\cF_i$ is completely separating.
        
    The above family of queries solves Model 3 but it cannot be used for Model 2, since the queries do not form a partition of the underlying set. Let us prove now the lower bound $r\binom{k}{\lfloor k/2\rfloor}$ in the case of Theorem \ref{thm2}, $r>2$ and $k>2$. Again, on each $Y_i$ Player $B$ picks a completely separating family $\cF_i=\{X_1^{i},\dots,X_k^{i}\}$, but now we specify that the dual of $\cF_i$ is the family of $\lceil k/2\rceil$-element sets of the $k$ element base set. Queries $Q_1, Q_2,\dots, Q_k$ obtained as above are still not partitions of $X$, so we have to modify them.  
    
    There are some elements of $Y_i$ that are not in $X_j^{i}$. Let us first add these elements to $X_j^\ell$ for some $\ell\neq i$. Repeating this for every $i$, we obtain a new query $Q'_j$ that is a partition of $X$ to $U_j^1,\dots, U_j^r$. Let us repeat this for every $j$. Then, as all elements of $U_j^{i}\setminus X^{i}_j$ are not in $Y_i$, the sets $U_1^{i},\dots, U_k^{i}$ still form a completely separating family restricted to $Y_i$. If Player $C_i$ would know that the defective element is in $Y_i$, then he could identify it. However, it is possible that for some $x\in Y_i$ and $y\not\in Y_i$, each set from $\{U_1^{i},\dots, U_k^{i}\}$ that contains $x$ also contains $y$. If that is the case and $x$ is the defective, then Player $C_i$ cannot rule out the possibility that $y$ is the defective.
    
    Observe that each element of $Y_i$ appears in $\lceil k/2\rceil$ sets $X_j^{i}$ by the choice of $\cF_i$. Therefore, such an element $x\in Y_i$ appears in  $\lceil k/2\rceil$ sets $U_j^{i}$. If $y\not\in Y_i$, then $y\in Y_\ell$ for some $\ell$, thus such an $y$ appears in $\lceil k/2\rceil$ sets $U_j^{\ell}$, thus $y$ appears in at most $\lfloor k/2\rfloor$ sets $U_j^{i}$. Thus, if $k$ is odd, then we cannot have that each set $U_j^{i}$ that contains $x$ also contains $y$ (as $\lfloor k/2\rfloor<\lceil k/2\rceil$), which completes the proof if $k$ is odd.

    
    If $k$ is even, then we still want that every element of $Y_{\ell}$ appears in less than $k/2$ sets $U_j^i$ for every $\ell\neq i$. To do this we extend the sets $X_j^{i}$ to $U_j^{i}$ more carefully. Then for each $j$ we take the elements of $Y_i\setminus X_j^{i}$, and add each of them to one of two different sets $U_j^{\ell}$ and  $U_j^{\ell'}$ in the following way. For each element $y\in Y_i$, we first consider the smallest $j$ such that $y\not\in X_j^{i}$, and then we add $y$ to $U_j^{\ell}$. Then we take the second smallest index $j'$ such that $y\not\in X_{j'}^{i}$, and then we add $y$ to $U_{j'}^{\ell'}$. Adding $y$ to these two sets alternately this way, any $y\not\in Y_i$ appears in less than $k/2$ sets $U_j^{i}$. Therefore, we cannot have that each set $U_j^{i}$ that contains $x$ also contains $y$ as every $x\in Y_i$ appears in $k/2$ sets $U_j^{i}$. We can execute this plan if $r\ge 3$ (to pick $i$, $\ell$ and $\ell'$) and $k>2$ (to pick $\ell$ and $\ell'$ among the $\lceil k/2\rceil$ sets not containing $y$).

    The case $r=2$ is dealt with in Proposition \ref{prop1}. 
    
    It is left to show that $g_2(2,r)\ge 2r-1$. Let $X=\{x_1,\dots,x_{2r-1}\}$.
    We define the Queries $Q_1=(X^1_1,X^1_2,\dots,X^1_r)$, $Q_2=(X^2_1,X^2_2,\dots,X^2_r)$ by $X^1_i=\{x_{2i-1}\}, X^2_i=\{x_{2i}\}$ for $i\le r-1$ and $X^1_r=X\setminus \cup_{i=1}^{r-1}X^1_i$, $X^2_r=X\setminus \cup_{i=1}^{r-1}X^2_i$. If the defective element is $x_j$ with $j<2r-1$, then it belongs as a singleton to a Player, thus that Player can identify it. If the defective element is $x_{2r-1}$, then Player $C_r$ can identify it, since $x_{2r-1}$ is the only element contained in both $X^1_r$ and $X^2_r$.
\end{proof}

In multiple proofs next, Players $B$ and $C$ can decide their strategy before the algorithm starts. They will apply modifications of the following strategy that we will call the \textit{Basic Strategy}. Let us identify elements of $X$ with 0-1 sequences of length $\lceil\log n\rceil$, and let $S_i$ denote the set of sequences whose $i$th bit is 1. In the Basic Strategy, the $i$th query of Player $B$ is $S_i$, and the $i$th query of Player $C$ is $\overline{S_{\lceil\log n\rceil -i+1}}$. In other words, Player $B$ asks whether the $i$th bit is 1, while Player $C$ asks whether the $i$th bit from the other direction is 0. 

The main advantage of a predetermined strategy is that whenever a YES answer is obtained to a query of $B$, then  $C$ obtains this answer, and also knows what the earlier queries of $B$ were. Since $C$ knows when the answers were YES, he knows the answer to all the queries of $B$ (and vice versa). The reason that $C$ asks the complements is that in this way when they ask about the same bit, one of them obtains the answer YES, thus one of them obtains the information.

Recall that Players $B$ and $C$ are not allowed to communicate during the algorithm. However, they can send messages to each other via asking certain queries. One simple way to send such a message is to ask $X$ as a query (which clearly would not serve any other purpose). Since the answer to this is YES, the other player will obtain this message in the models studied. One can ask the same query multiple times, thus any message could be sent. 
These messages will be added at some points during the Basic Strategy.

We continue with the proof of Proposition \ref{prop5}. Recall that it deals with the model where Player $A$ decides who asks the next query, the Questioner obtains the answer, while the other Player obtains the answer (and the query) if the answer is YES. The statement is that at least one player can identify the defective element after $\log n+1$ queries.

\begin{proof}[Proof of Proposition \ref{prop5}]
    Players $B$ and $C$ apply the Basic Strategy. Then at query $\log n+1$ they arrive to the same bit $i$. Let us assume without loss of generality that the defective element is in $S_i$, thus the answer to the query by Player $B$ is YES. Both players know this answer. No matter which one of them arrived to the last query, that player knows the answer to every earlier query, thus can identify the defective element. Indeed, he knows the answer to his own queries by the properties of Model 5, and knows the answer to the queries asked by the other player, since he knows that the other player asked the queries till the $i$th bit, and knows which of those queries were answered YES.
\end{proof}

The inequality in this Proposition is not sharp. We determined the exact number of question for values of $n$ up to 20:

\begin{equation*}
    f_5(n)=
    \begin{cases}
        1 & \text{if $n = 2$} \\
        2 & \text{if $n = 3$}\\
        3 & \text{if $4 \leq n \leq 6$} \\
        4 & \text{if $7 \leq n \leq 11$} \\
        5 & \text{if $12 \leq n \leq 20$.} \\
    \end{cases}       
\end{equation*}


Let us continue with the proof of Theorem \ref{thm6}. Recall that  Model 6 differs from Model 5 in that both players should be able to identify the defective element. The bound we have to prove is  $f_6(n)\le \lceil\log n\rceil+2\lceil\sqrt{\log n}\rceil+2$.

\begin{proof}[Proof of Theorem \ref{thm6}]
    Let $P$ be the set of numbers of the form $i\lceil \sqrt{\log n}\rceil$ for $i\le \lceil \sqrt{\log n}\rceil$. Players $B$ and $C$ follow the Basic Strategy with some additional queries. For simplicity, we describe the strategy for Player $B$; it is analogous for Player $C$. Whenever Player $B$ asks $S_p$ with $p\in P$, the next query of Player $B$ is a \textit{special query}. Afterwards Player $B$ continues with the Basic Strategy, i.e., her next query is $S_{p+1}$. The special query is $X$, thus the answer is always YES, thus Player $C$ obtains this information. Therefore, Player $C$ also knows which special query this was, thus Player $C$ knows the answer to each $S_i$ with $i<p$.

    Whenever a query corresponding to the $i$th bit is asked by both players, one of them obtains a YES answer, thus the other player knows every bit as in the proof of Proposition \ref{prop5}. Assume now that the special query at $p\in P$ is asked by both players. Then Player $B$ knows the answer to each $S_i$ with $i\le p$ since he asked those queries, and the answer to each $\overline{S_i}$ with $i\ge p$ since those were asked by Player $C$ and this special query shows this fact. Player $C$ analogously knows each bit, thus the algorithm is finished.

    As there are $\lceil \sqrt{\log n}\rceil$ elements in $P$, there are at most $\lceil \sqrt{\log n}\rceil+1$ special queries asked. For each $i\le \log n$, one of the queries $S_i$ and $\overline{S_i}$ is asked, except for the the values $i$ between two consecutive elements of $P$. Indeed, if without loss of generality Player $B$ asks the special query at $p$ first, it is possible that he asks all the queries between $p$ and the next element $p+\lceil \sqrt{\log n}\rceil$ of $P$ before Player $C$ asks the special query at $p$. However, if Player $B$ would ask also the special query at $p+\lceil \sqrt{\log n}\rceil$, then this special query would be asked by both players before the special query at $p$, thus the algorithm would finish before Player $C$ would ask the special query at $p$, a contradiction. Therefore, there are at most $\lceil \sqrt{\log n}\rceil+1$ entries $i$ with both $S_i$ and $\overline{S_i}$ asked, hence there are  $\lceil\log n\rceil+\lceil\sqrt{\log n}\rceil+1$ non-special queries asked.    
\end{proof}

We remark that we could get rid of the $+2$ term  by making the strategy less symmetric: Player $B$ would ask the special query at $p$ after asking $S_p$ and Player $C$ would ask it before asking $\overline{S_p}$.

Let us continue with the proof of Theorem \ref{thm4}. Recall that the player asking the query does not get any information, but it is enough if one of Players $B$ and $C$ finds the defective element. The theorem states that $f_4(n)\le \lceil\log n\rceil+\lceil\log\log n\rceil+2$.


\begin{proof}[Proof of Theorem \ref{thm4}]
Players $B$ and $C$ again follow the Basic Strategy with the following alterations. Throughout the process, they maintain a left endpoint $\ell$, a right end point $r$, and a target position $t=\frac{\ell+r}{2}$.
So at the beginning $\ell=1$ and $r=\lceil\log n\rceil$. $B$ plays along the Basic Strategy, unless he is to ask a query $S_i$ with $i>t$, while $C$ plays along the Basic Strategy unless he is to ask a query $\overline{S_i}$ with $i<t$. At these moments $B$ and $C$ ask a \textit{special query} for which the answer is always YES and which is of the form $\cap_{i\in I}S_i\cap \cap_{i\in I'}\overline{S_i}$. The special query  is the smallest set $S$ the questioner knows the defective element belongs to. Formally, without loss of generality assume it is $C$ who asks the special query $Q=Q_1\cap Q_2$, where $S_1$ is the set $C$ knows because of previous queries of $B$, and $S_2$ is the set $C$ knows about his own queries. More precisely, if we denote by $YES_B$ the set of queries asked by $B$ that were answered YES, then
$$Q_1=\bigcap_{S'\in YES_B}S'\cap\bigcap_{\exists i'>i: S_{i'}\in YES_B}\overline{S_i}.$$
Observe that $C$ does know $\bigcap_{S'\in YES_B}S'$ as these are queries of $B$ that were answered YES, so he received them. He also knows $\bigcap_{\exists i'>i: S_{i'}\in YES_B}\overline{S_i}$ as the players agreed to play along the Basic Strategy, so if $C$ did not receive $S_{i}$ but received an $S_{i'}$ with $i<i'$, then $C$ knows that the defective element is in $\overline{S_i}$. An important consequence of asking the special query $Q$ that has a part $Q_1$ is that $B$ after receiving $Q$ will know the answers to all of his previous queries. Indeed, as $C$ asks $\bigcap_{S'\in YES_B}S'$, $B$ knows which of his queries have been answered YES, and thus all other queries must have been answered NO. If before $C$ is asking a special query, there had been special queries from $B$, then $C$ knows the answers to those of his queries that were asked before $B$'s last special query, so $C$ can include these $S_j$s and $\overline{S_j}$s into $Q_2$.

There are two ways how a special query can be asked. The simpler case is if $B$ asks this instead of asking $S_{\lfloor t\rfloor+1}$ (or similarly $C$ asks this instead of asking $S_{\lceil t\rceil-1}$). In this case, $\ell$ is changed to $\lfloor t\rfloor+1$ and $t$ is changed according to the rule $t=\frac{\ell+r}{2}$. (If $C$ was asking the special query, then $r$ becomes $\lceil t\rceil-1$ and $t$ is updated accordingly.) Another possibility is that a Player asked a special query and that moved the target $t$ ``behind" the other player. For example if $B$ is first to reach the bit $\frac{1}{2}\log n$, but by that time $C$ has already asked $\overline{S_{\frac{5}{8}\log n}}$, then $t$ becomes $\frac{3}{4}\log n>\frac{5}{8}\log n$. This situation can only occur when say, $C$ asks his first query after a special query of $B$ (maybe $B$ has asked queries since, but not $C$), so as mentioned above, at this moment $C$ knows the answers to all the previous queries by himself. Therefore, the $Q_2$ part of his special query will reveal his current position $i$ to $B$, so the two Players can set $r$ to be $i$ and then update $t$.

As the value of $r-\ell$ obtains at least halved after each change to $t$, we obtain that the number of special queries is at most $\lceil\log (\lceil \log n\rceil)\rceil$. To finish the prof, we need to analyze the situation when $B$ has asked $S_i$, $C$ has asked $\overline{S_{i+1}}$, so until this point $\lceil \log n\rceil$ normal and at most $\lceil\log (\lceil \log n\rceil)\rceil$ special queries have been asked.

\medskip

\textsc{Case I.} We have $t>i+1$, i.e., that target is behind $C$.

\smallskip

Then by the above observation, $C$ knows the answers to all his previous queries. If he is to ask the next query (a special query as $t>i+1$), then $B$ will learn the answers to all bits and thus will know the defective element. If $B$ is to ask $S_{i+1}$, then there are two possibilities. If the answer is YES, then $C$ will receive this, and will know the defective element, as because of the Basic strategy, $C$ will know all the YES and thus all NO queries of $B$. Finally, if the answer to $S_{i+1}$ would be NO, then previously the answer to $C$'s query $\overline{S_{i+1}}$ must have been YES and thus received by $B$. So $B$ instead of asking $S_{i+1}$ can ask a special query ensuring that $C$ will know the defective element.

\smallskip

\textsc{Case II.} We have $t=i+\frac{1}{2}$ or $t=i+1$.

\smallskip

If $C$ is to ask the next query, then the query is a special one that lets $B$ know the answers to all his queries and moves the target behind $B$. So any later query from $B$ would be a special query letting $C$ know the defective element. If after the special query of $C$, $C$ is to ask the next query, then that query is $\overline {S_i}$. Then the situation is analogous to the second subcase of Case I: if the answer is YES, then $B$ receives this and finds out the answer to the all the queries by $C$, while if the answer is NO, then $C$ already knew that since $B$ had asked $S_i$, so he can ask a special query instead of $\overline {S_i}$.

If $B$ is to ask the next query and $t=i+\frac{1}{2}$, then the roles of $B$ and $C$ are symmetric, so the above paragraph can be applied. Finally, if $B$ is to ask the next query and $t=i+1$, then after the next query both players would proceed with a special query and we finish as before.

In all cases, we finished with at most 2 extra queries. This completes the proof.
\end{proof}

Let us continue with the proof of Theorem  \ref{thm7}. Recall that now Players $B$ and $C$ can agree in a strategy beforehand, but the elements are indistinguishable at that point. They alternate asking queries, and at least one of them should be able to identify the defective element.

\begin{proof}[Proof of Theorem \ref{thm7}]
    Recall first that for the golden ratio we have that $1/\phi+1/\phi^2=1$.
    
    Let us start with the upper bound. At each point, one of Players $B$ and $C$ is the \textit{main questioner}, the other is the \textit{auxiliary questioner}. Let us assume that Player $B$ asks the first query, then he starts as the main questioner. They change roles whenever there is a YES answer.

    If the answer to the query of the main Questioner is NO, the auxiliary questioner asks the whole underlying set $X$. This tells the main Questioner that the answer to his Query was NO.

    The main questioner always maintains the set of possibly defective elements $X_i$ and asks a subset $Y_i$ of $X_i$ of order $\lfloor |X_i|/\phi\rfloor$ as a query. If the answer is YES, the other player becomes the main questioner, and $X_{i+1}:=Y_i$ is the set of possibly defective elements. If the answer is NO, then the auxiliary questioner queries $X$, and the main Questioner sets the set of possibly defective elements to $X_i\setminus Y_i$. This way, (omitting the floor signs) the size of the set of possibly defective elements either decreases to $|X_i|/\phi$ with one query or decreases to $|X_i|(1-1/\phi)=|X_i|/\phi^2$ with two queries. Clearly this implies the upper bound $f_7(n)\le \log_\phi n$.

    Let us continue with the lower bound. We describe the strategy of the adversary recursively. Note that usually an adversary strategy means that we pick the answers so that there is still a possible solution. However, as the elements are indistinguishable, we can also control the intersection of the queries by the two players in the case the second player does not have any information on the query of the first player. For example, if the first player asks a set of size $n/2$ and the answer is NO, and the second player asks a set of size $n/2$, the adversary can decide that they asked the same set.
    
    If the first query has size at least $\lfloor n/\phi\rfloor$, the answer is YES, otherwise NO. If the first answer is YES, he applies the same strategy starting at the next query. If the first answer is NO, then note that the complement of the first query has size at least $\lfloor n/\phi^2\rfloor$. If the second query has size at least $\lfloor n/\phi^2\rfloor$, then the adversary decides that the intersection of the second query with the complement of the first query has size at least $\lfloor n/\phi^2\rfloor$ and the answer is YES. Moreover, he tells both players all the answers to the queries. This way after two queries there are at least $\lfloor n/\phi^2\rfloor$ possibly defective elements, and the adversary applies the same strategy starting at the next query. On the other hand, if the second query has size less than $\lfloor n/\phi^2\rfloor$, the adversary answers NO and decides that the union of the two queries has size at most $\lfloor n/\phi\rfloor$ and tells this information to both players. Again, after two queries there are at least $\lfloor n/\phi^2\rfloor$ possibly defective elements, and the adversary applies the same strategy starting at the next query. This completes the proof.
\end{proof}

\section{Concluding remarks}

There are other models that one can obtain by combining the ideas used in the  paper. Some of these do not make much sense, e.g. if $A$ decides who asks the next query, and the player asking the query does not get any information, we cannot require that both players identify the defective, since $A$ can make sure that only one player asks all the queries, and that player can obviously never identify the defective. 

Other models are just not of much interest, because their solution is trivial. E.g., one might ask what happens in Model 1 if we require that both players can identify the defective. This cannot be easier than the simple puzzle we started with, where Player $C$ has to identify the defective. But the solution to that is a completely separating family, and it is not hard to see that if the queries form a completely separating family, then Player $D$ can also identify the defective element, completely solving this case. Another example is if Players $B$ and $C$ ask queries alternately, the other player obtains the YES answers and the elements are distinguishable. Then the players agree on a strategy and each player knows the answer to the queries of the other player, thus they can always decrease the set of possibly defective elements by a factor of 2, which leads to an algorithm requiring $\log n$ queries.

We never forbade agreeing in a strategy before the start of an adaptive algorithm. The reason is that forbidding this leads to a sort of philosophical question. Can the players agree in a strategy in constant time? They can communicate during the algorithm (by sending messages via asking certain queries) and send any number $k$ to each other, e.g., by asking the query $X$ exactly $k$ times in a row. Without discussing the strategy, but assuming optimal play by the players, can we assume that the other player understands that the message is $k$ and can he decode this message? If yes, then forbidding the discussion beforehand only gives a constant additive term.

\smallskip

In this paper we restricted our study to worst case bounds. Yet one could also consider the expected length of algorithms in case the defective element is chosen randomly (and then we can also decide if Player $A$ should choose randomly which Player can ask a query next in the appropriate models, or this could be the choice of $A$). 

\smallskip

Let us finally remark that while cooperation was essential in our algorithms, we could have studied competitive versions of Combinatorial Search instead, using the same models. For example, in Models 4,5,6,7,8 it could be the goal of Player $B$ that Player $C$ should not be able to identify the defective element (and vice versa). What are the chances that Player $B$ wins when the defective is a (uniformly) random element and how many queries are needed in the optimal algorithm? In Models 1,2,3 the goal of Player $B$ could be e.g. that one player cannot identify the defective element but multiple players can.

\bigskip

\textbf{Funding}: Research of Gerbner is supported by the National Research, Development and Innovation Office - NKFIH under the grants SNN 129364, FK 132060, and KKP-133819.

Research of Keszegh is supported by the J\'anos Bolyai Research Scholarship of the Hungarian Academy of Sciences, by the National Research, Development and Innovation Office -- NKFIH under the grant K 132696 and FK 132060, by the \'UNKP-21-5 and \'UNKP-22-5 New National Excellence Program of the Ministry for Innovation and Technology from the source of the National Research, Development and Innovation Fund and by the ERC Advanced Grant “ERMiD".

Research of Patk\'os is supported by the National Research, Development and Innovation Office - NKFIH under the grants SNN 129364 and FK 132060.

 This research has been implemented with the support provided by the Ministry of Innovation and Technology of Hungary from the National Research, Development and Innovation Fund, financed under the  ELTE TKP 2021-NKTA-62 and the TKP2021 (project no. BME-NVA-02) funding scheme.

\end{document}